\input amstex
\documentstyle{amsppt}
\font\thstil=cmsl10

\def\ltextindent#1{\hbox to x
\hangindent{#1\hss}\ignorespaces}

\newskip\proofskipamount
\proofskipamount=8pt plus 2pt minus 2 pt

\def\endth{\par\ifdim\lastskip<\bigskipamount
\removelastskip\penalty55
\bigskip\fi}

\def\ifundefined#1{\expandafter\ifx
\csname#1\endcsname\relax}

\newif\ifdevelop \developfalse

\newtoks\chnumber
\newtoks\sectionnumber
\newcount\equationnumber
\newcount\thnumber

\def\assignnumber#1#2{%
\ifundefined{#1}\relax\else\message{#1 already defined}\fi
\expandafter\xdef\csname#1\endcsname
{\if-\the\sectionnumber\else\the\sectionnumber.\fi\the#2}%
}%

\def\beginsektion #1 #2 {\vskip0pt plus.1\vsize\penalty-250
\vskip0pt plus-.1\vsize
\bigbreak\bigskip
\sectionnumber{#1} \equationnumber0\thnumber0
\noindent{\bf #1. #2}\par
\nobreak\medskip\noindent}

\def\nologo{\expandafter\let
\csname logo\string @\endcsname=\empty}
\def\:{:\allowbreak }
\def\eq#1{\relax
\global\advance\equationnumber by 1
\assignnumber{EN#1}\equationnumber
{\tenrm (\csname EN#1\endcsname)}
}

\def\eqtag#1{\relax\ifundefined{EN#1}
\message{EN#1 undefined}{\sl (#1)}%
\else(\csname EN#1\endcsname)\fi%
}

\def\thname#1{\relax
\global\advance\thnumber by 1
\assignnumber{TH#1}\thnumber
\csname TH#1\endcsname
}

\def\beginth#1 #2 {\bigbreak\noindent{\bf #1\enspace \thname{#2}}
    ---\hskip4pt}
\def\thtag#1{\relax\ifundefined{TH#1}
\message{TH#1 undefined}{\sl #1}%
  \else[\csname TH#1\endcsname]\fi}

\def\si{\sigma}

\def\AA{{\Bbb A}}

\def\CC{{\Bbb C}}

\def\PP{{\Bbb P}}

\def\RR{{\Bbb R}}
\def\TT{{\Bbb T}}
\def\ZZ{{\Bbb Z}}

\def\FSA{\Cal A}
\def\FSB{\Cal B}

\def\FSD{\Cal D}

\def\FSH{\Cal H}

\def\FSU{\Cal U}

\def\FSZ{\Cal Z}

\def\goa{{\goth a}}
\def\goh{{\goth h}}
\def\gog{{\goth g}}
\def\gok{{\goth k}}
\def\gol{{\goth l}}

\def\gop{{\goth p}}
\def\gos{{\goth s}}

\def\goS{{\goth S}}

\def\id{\operatorname{id}}
\def\phi{\varphi}

\def\Aq{\FSA_q}
\def\Uq{\FSU_q}
\def\Bq{\FSB_q}
\def\detq{\det\nolimits_q}
\def\Lp{L^+}
\def\Lm{L^-}

\def\Hst{\FSH^{(\si,\tau)}}
\def\xst{x^{(\si,\tau)}}
\def\AT{\FSA(\TT)}
\def\Uh{\Uq(\goh)}

\def\End{\operatorname{End}}
\def\Hom{\operatorname{Hom}}

\def\innp#1 #2 {\langle #1 , #2 \rangle}
\let\eps=\varepsilon

\let\ten=\otimes
\def\AA{AA}
\def\AW{AW}
\def\Chin{Ch}
\def\DKa{DK1}
\def\DKb{DK2}
\def\DN{DN}
\def\DNS{NDS}
\def\Dr{Dr}
\def\GR{GR}
\def\Hay{Ha}
\def\HS{HS}
\def\Hela{H1}
\def\Helb{H2}
\def\Jima{J1}
\def\Jimb{J2}
\def\Jimc{J3}
\def\KWa{K1}
\def\KWb{K2}
\def\KWc{K3}
\def\KWd{K4}
\def\VK{VK}
\def\Ku{Ku}
\def\Lz{Lz}
\def\Ma{M1}
\def\Mb{M2}
\def\Mc{M3}
\def\Mas{Mas}
\def\Nmac{N}
\def\NSa{NS1}
\def\NSb{NS2}
\def\NYM{NYM}
\def\Po{Po}
\def\FRT{RTF}
\def\Ro{R}
\def\St{St}
\def\Su{Su}
\def\UT{UT}
\def\VS{VS}
\def\Wo{Wo}
\topmatter
\title Some remarks on the construction of quantum symmetric spaces
\endtitle
\author Mathijs S.\ Dijkhuizen
\endauthor
\affil Department of Mathematics, Faculty of Science,  \\
Kobe University, Rokko, Kobe 657, Japan
\endaffil
\email msdz@math.s.kobe-u.ac.jp
\endemail
\thanks The author acknowledges financial support
by the Japan Society for the Promotion of Science (JSPS)
and the Netherlands Organization for Scientific Research (NWO).
\endthanks
\abstract
We present a general survey of some recent  developments 
regarding the construction of compact 
quantum symmetric spaces and the analysis
of their zonal spherical functions in terms of $q$-orthogonal
polynomials. In particular, we define a one-parameter family
of two-sided coideals in $\Uq(\gog\gol(n,\CC))$ and express the
zonal spherical functions on the corresponding quantum projective
spaces as Askey-Wilson polynomials containing two continuous
and one discrete parameter.
\endabstract
\date 31 May 1995 \enddate
\rightheadtext{Quantum symmetric spaces}
\leftheadtext{Quantum symmetric spaces}
\endtopmatter
\document
\beginsektion 0 {Introduction}
In this paper, we discuss some recent progress in the analysis
of compact quantum symmetric spaces and their zonal spherical functions.
It is our aim to present a more or less coherent overview of a number of
quantum analogues of symmetric spaces that have recently been 
studied. We shall try to emphasize those  aspects 
of the theory that distinguish quantum symmetric spaces from
their classical counterparts.

It was recognized not long after the introduction of
quantum groups by Drinfeld \cite{\Dr}, Jimbo \cite{\Jima}, and
Woronowicz \cite{\Wo}, that the quantization of symmetric
spaces was far from straightforward. We discuss some of the main
problems that came up.

First of all, there was the, at first sight  rather annoying, lack
of interesting quantum subgroups. Let us mention, for instance, the
well-known fact that there is no analogue of $SO(n)$ inside the
quantum group $SU_q(n)$. A major step forward in overcoming this
hurdle was made by Koornwinder \cite{\KWb}, who was able to construct
a quantum analogue of the classical 2-sphere $SU(2)/SO(2)$ and express 
the zonal spherical functions as a subclass of Askey-Wilson polynomials 
by
exploiting the notion of a twisted primive element in the quantized 
universal enveloping algebra $\Uq(\gos\gol(2,\CC))$. This 
``infinitesimal'' approach was considerably generalized by
Noumi in his fundamental paper \cite{\Nmac}, where he used the more
general notion of a two-sided coideal in $\Uq(\gog)$ to study
quantum analogues of the symmetric spaces $SU(n)/SO(n)$ and
$SU(2n)/Sp(n)$ and express their zonal spherical functions
as Macdonald's symmetric polynomials associated with root system
$A_{n-1}$ (cf.\ \cite{\Ma}). In the same paper, 
both the $L$-operators
(cf.\ \cite{\FRT}) and constant solutions of the
reflection equation (cf.\ \cite{\Ku}) 
were used for the first time to construct
suitable two-sided coideals.

A second rather striking phenomenon was the appearance of
parametrized families of quantum symmetric spaces with one
and the same classical counterpart. Podle\' s \cite{\Po}
was the first to exhibit a continuously parametrized 
family of quantum analogues of the classical 2-sphere.
The zonal spherical functions on these quantum
spheres were completely analysed by 
Koornwinder \cite{\KWb}, \cite{\KWd}, who used the infinitesimal
method.

In this paper, we announce some recent results by M.\ Noumi and
the author showing that the parameter phenomenon is not
restricted to 2-spheres but extends to complex
projective spaces of arbitrary dimension and even beyond.
The zonal spherical functions on our 
family of quantum projective spaces are expressed
as Askey-Wilson polynomials containing two continuous and
one discrete parameter.

A third difference 
is the rather {\sl ad hoc} approach to quantum
symmetric spaces. As far as the author knows, no method has been
devised so far that provides a unified approach to even a subclass
of quantum symmetric spaces. Related to this are the computational
difficulties involved in the quantum case,
especially in the computation of the
radial part of a quantum Casimir operator. First applied
by Koornwinder \cite{\KWd} in the rank one case, then used by Noumi
\cite{\Nmac} in the general rank case (see also \cite{\UT}), this method
has by now become one of the standard means to identify zonal
spherical functions on quantum symmetric spaces.

The organization of this paper is as follows. In section 
1 we briefly review the theory of classical compact
symmetric spaces and their zonal spherical functions. In sections
2 and 3 we recall some basic facts concerning
compact quantum groups, particularly the quantum unitary group,
and sketch one possible general method to analyse quantum
symmetric spaces. In section 4 we announce some
recent results by M.\ Noumi and the author on a family
of quantum projective spaces. Finally, in section
5 we briefly discuss some patterns
underlying the study of higher rank quantum symmetric spaces.

The author would like to express his
sincere gratitude to Prof.\ M.\ Noumi and T.\ Sugitani.
Many of the general ideas expressed in this paper originated
with them, still others ripened in the countless discussions
we had during our informal seminars in Tokyo. The
author would also like to thank Prof. Tom H. Koornwinder,
who introduced him to the subject of harmonic analysis on
quantum groups.
\beginsektion 1 {Classical compact symmetric spaces}
The material treated in this section is well-known by now
(cf.\ Helgason \cite{\Hela}, \cite{\Helb}, Heckman and Schlichtkrull
\cite{\HS}).

Let $G$ be a compact connected simple Lie group.
To simplify our statements we assume that $G$ is simply
connected, but this is not essential.
 Let $K\subset
G$ be a closed subgroup such that there exists a (necessarily
unique) involutive
automorphism $\theta\colon G\to G$ with $G^\theta_0 \subset
K\subset G^\theta$. Here $G^\theta$ denotes the subgroup of
fixed points of $\theta$, $G^\theta_0$ its
connected component of the identity. Then the pair $(G,K)$ is
an irreducible  Riemannian symmetric pair, and the $G$-homogeneous space
$G/K$ is an irreducible Riemannian symmetric space of compact type I 
(in the terminology of Helgason \cite{\Hela}).

For any (finite-dimensional continuous) irreducible representation $\pi$
of $G$,  the trivial representation $\delta$ of $K$ occurs at
most once in the irreducible decomposition of the 
restriction of $\pi$ to $K$. In other words, $(G,K)$ is a 
{\sl Gelfand pair}.
 Those irreducible representations
$\pi$ of $G$ that have non-zero $K$-fixed vectors are
called {\sl spherical representations}. 

Let $\FSD(G/K)$ denote the algebra of left $G$-invariant
differential operators on the space $G/K$. A {\sl (zonal)
spherical function} is by definition a $K$-biinvariant
differentiable function on $G$ which is a joint eigenfunction of
all differential operators $D\in \FSD(G/K)$. 

There is a 1-1 correspondence between spherical representations
and spherical functions (determined up to a scalar multiple).
 To make this more explicit, we write $\gog\subset \gog_\CC$ for
the Lie algebra of $G$ and its complexification. Let $\goh_\CC\subset
\gog_\CC$ be a Cartan subalgebra, $R:= R(\gog_\CC, \goh_\CC)$
the corresponding root system. We fix a choice of positive
roots $R^+\subset R$. Let $P^+\subset P\subset \goh_\CC^\ast$
denote the corresponding cone $P^+$ of dominant weights in
the weight lattice $P$. The irreducible representation of $G$ with highest
weight $\lambda\in P^+$ is denoted by $V(\lambda)$.
 Let $P_K^+\subset P^+$ denote the subset of
dominant weights 
corresponding to the spherical representations.
If we write $\FSH$ for the algebra of $K$-biinvariant
(continuous) representative (i.e.\ $G$-finite) functions
on $G$, then there is the canonical decomposition
$$
\FSH = \bigoplus_{\lambda \in P^+_K} \FSH(\lambda),
\eqno\eq{PclHdecomp}
$$
where $\FSH(\lambda)$ is the intersection of $\FSH$ and the subspace
$W(\lambda)\subset L^2(G)$ spanned by the matrix 
coefficients of the representation $V(\lambda)$.
Each of the subspaces $\FSH(\lambda)$ is one-dimensional.
If we choose a non-zero $\phi_\lambda\in \FSH(\lambda)$,
then $\phi_\lambda$ will be a spherical function. The right 
$G$-translates of $\phi_\lambda$ span a subspace which is equivalent
as a representation of $G$ to $V(\lambda)$.
All spherical functions can be
obtained this way. 

In order to give a more explicit description of the spherical functions
$\phi_\lambda$, we recall the definition of (generalized) Jacobi polynomials
(cf.\ \cite{\HS}).
Let $\Sigma\subset V$ be a (possibly non-reduced) root system of rank $l$
in a real vector space $V$.
We fix an inner product $\langle \cdot\,,\, \cdot \rangle$ on $V$ which is
invariant under the Weyl group $W=W_\Sigma$.
Via this inner product we may identify $V$ with its dual $V^\ast$.
We fix a choice $\Sigma^+\subset \Sigma$ of positive roots. We then have
 a cone of dominant weights 
$P^+_\Sigma$ inside the weight lattice $P_\Sigma$.
 The symbol $\preceq$ denotes the usual
dominance ordering of weights. Let $\CC[P_\Sigma]$ 
denote the group algebra of
the free abelian group $P_\Sigma$ with canonical elements
$e^\lambda\in \CC[P_\Sigma]$ ($\lambda\in P_\Sigma$). The Weyl group
 $W$ acts on $\CC[P_\Sigma]$ by
$w\cdot e^\lambda = e^{w\lambda}$. The orbit sums
$m_\lambda := \sum_{\mu\in W\lambda} e^\mu$ ($\lambda \in P^+_\Sigma)$ are
$W$-invariant and form a basis of the subalgebra
$\CC[P_\Sigma]^W$ of Weyl group invariants. Let $Q^\vee\subset V^\ast$
denote the dual root lattice. Then $T:= V^\ast/2\pi Q^\vee$ is
a compact torus. If we interpret $e^\lambda\in\CC[P_\Sigma]$ as a function
on $T$ by $e^\lambda(x) := e^{i\langle \lambda\,,\, x \rangle}$, then 
the group algebra $\CC[P_\Sigma]$ is identified with the algebra $\FSA(T)$
of representative functions on the torus $T$. This identification 
is compatible with the natural actions of $W$ on both algebras.

Let now $k:\Sigma \to [0,\infty), \, \alpha\mapsto k_\alpha$ 
be a $W$-invariant (so-called {\sl multiplicity})
function. 
We define a continuous weight function on $T$ and the corresponding
inner product on functions of $T$:
$$
\delta_k(x):= \prod_{\alpha\in \Sigma^+} |2\sin({1\over 2}\langle \alpha,
x\rangle)|^{2k_\alpha}, \quad
\langle f, g \rangle_k := \int_T f(x) \overline{g(x)} \delta_k(x) 
\hbox{d}x,
\eqno\eq{clweight}
$$
where $\hbox{d}x$ is the normalized Lebesgue measure on $T$.
The {\sl Jacobi polynomial} $P_\lambda^k\in \CC[P_\Sigma]^W$ 
of degree $\lambda
\in P^+_\Sigma$ is defined as the 
element of the form $\sum_{\mu\preceq\lambda}
c_{\lambda,\mu} m_\mu$ with $c_{\lambda,\lambda}=1$ such that
$$
\langle P^k_\lambda, m_\mu\rangle_k = 0, \quad \mu\in P^+_\Sigma, 
\mu\prec\lambda.
$$
This last condition is equivalent to the requirement that
$P_\lambda^k$ satisfies the following second-order partial
differential equation:
$$
\left ( \Delta + \sum_{\alpha\in \Sigma^+} k_\alpha \cot({1\over 2}
\langle \alpha, x\rangle)\partial_\alpha\right ) P^k_\lambda (x)
= -\langle \lambda, \lambda + \sum_{\alpha\in \Sigma^+}
k_\alpha\, \alpha \rangle  P^k_\lambda(x).
\eqno\eq{jacdiff}
$$
Here $\Delta$ is the usual Laplace operator on the Euclidean
space $V$, $\partial_\alpha$ the partial derivative
in the direction of $\alpha$. We denote the partial differential operator
on the left-hand side of \eqtag{jacdiff} by $L(k)$.
Clearly, the Jacobi polynomials $P^k_\lambda$ ($\lambda\in P^+_\Sigma$)
form a basis of the algebra $\CC[P_\Sigma]^W$.
One can prove that 
$$
\langle P_\lambda,P_\mu\rangle = 0, \quad
\lambda, \mu\in P^+_\Sigma,\; \lambda\neq\mu.
$$
 Hence, the $P^k_\lambda$ are a family of orthogonal
polynomials. One can actually show that $L(k)$ is contained
in a certain commutative algebra $\FSD(k)$ of self-adjoint 
differential operators
generated by  $l$ algebraically independent generators.
The $P^k_\lambda$ can be characterized as the joint eigenfunctions
in $\CC[P_\Sigma]^W$ of the operators in $\FSD(k)$.

Let us now return to the setting of the symmetric space $G/K$. 
Let $\gok\subset \gog$ denote the Lie algebra of $K$.
The eigenspace decomposition 
of $\theta\colon \gog\to\gog$ is written as 
$\gog=\gok\oplus i\gop \subset \gog_\CC$. 
Let $\goa\subset\gop$ be a maximal
abelian subspace. Write $\Sigma'\subset \goa^\ast$
for the  restricted root system of $G/K$, 
$m_\alpha$ for the multiplicity
of $\alpha\in\Sigma$.

Let us put 
$$
\Sigma := 2\Sigma', \quad \quad k_{2\alpha} := {1\over 2}m_\alpha.
\eqno\eq{multdef}
$$
Then $k\colon \alpha \mapsto k_\alpha$ is a $W$-invariant
multiplicity function on the root system $\Sigma\subset \goa^\ast$.

Now $T:= \exp(i\goa)/\exp(i\goa)\cap K$ is a compact torus in $G$ 
isomorphic with the quotient $\goa/2\pi Q^\vee_\Sigma$ via the exponential
mapping. Restriction defines an algebra
isomorphism of $\FSH$ onto the algebra $\FSA(T)^W$ of $W$-invariant
representative functions on the torus $T$. 

Let $D\in \FSD(G/K)$ be a $G$-invariant differential operator
on $G/K$. It is obvious that $D$ maps $\FSH$ into itself.
There is a uniquely determined differential operator $\hbox{rad}(D)$
acting on functions on $T$ such that $\hbox{rad}(D)(f_{|T}) =
D(f)$. This differential operator is called the {\sl radial part}
of $D$. The mapping $\hbox{rad}$ actually induces an algebra
isomorphism of $\FSD(G/K)$ onto the algebra $\FSD(k)$, where 
the multiplicity function $k$ is defined in \eqtag{multdef}. 
Under this isomorphism,
the Laplace-Beltrami operator $\Delta$ on the Riemannian symmetric
space $G/K$ is mapped onto the operator $L(k)$.  Let  us recall
that $\Delta$ arises as the image of the standard Casimir  operator
on $G$ which can be viewed as a central element of
the universal enveloping algebra $\FSU(\gog_\CC)$.

Suppose now that the Cartan subalgebra $\goh_\CC\subset \gog_\CC$
is invariant under the complex-linear extension of $\theta$ to $\gog_\CC$,
and that, in addition, $\goa\subset \goh_\CC$.
Then we have the decomposition 
$\goh_\CC = \goa_\CC \oplus (\gok_\CC\cap \goh_\CC)$.
The representation $V(\lambda)$ ($\lambda\in P^+$) is
spherical if and only if $\lambda$ vanishes on $\gok_\CC
\cap \goh_\CC$ and the restriction $\mu$
of $\lambda$ to $\goa$ is a dominant weight in $P^+_\Sigma$.
The spherical function
$\phi_\lambda$ ($\lambda\in P^+_K$), when restricted to
$T$, is a scalar multiple of the Jacobi polynomial
$P^k_{\mu}$, where $k$ is defined as in \eqtag{multdef}.

Since the spherical functions $\phi_\lambda$ are matrix coefficients
of irreducible representations, the Schur orthogonality relations
imply that the $\phi(\lambda)$ are mutually orthogonal with respect
to the inner product on $L^2(G)$.
 It can be shown that this inner product 
coincides  on $\FSH$
with the inner product \eqtag{clweight} up to a scalar multiple.
\beginsektion 2 {Compact quantum groups}
In the remainder of this paper, we assume that $G$ is one of
the following compact Lie groups: $U(n)$, $SU(n)$ ($n\geq 2$),
$SO(n)$ ($n\geq 3$), $Sp(n)$ ($n\geq 2$). Actually, to simplify our
statements, we take 
$G=U(n)$ in this section.
 {\sl Mutatis mutandis} all the definitions and 
statements in this section are valid for the cases $SU(n)$, $SO(n)$, and
$Sp(n)$ too.

Let us fix  $0<q<1$ and $n\geq 2$. 
We identify $\gog_\CC$ with the Lie algebra $\gog\gol(n,\CC)$
of complex $n\times n$ matrices. As Cartan subalgebra 
$\goh_\CC\subset \gog_\CC$ we take the subspace of diagonal
matrices. Then, with the usual choice of $\eps_i\in \goh_\CC^\ast$, 
we can
identify the weight lattice $P$ with 
the free $\ZZ$-span of the $\eps_i$. 
There is a unique $\ZZ$-valued symmetric 
bilinear form $\langle \cdot \, , \,\cdot\rangle$ on
 $P$ such that
 $\langle \eps_i , \eps_j\rangle = \delta_{ij}$. 
Via this pairing, we shall identify $P$ 
with its dual $P^\ast = \Hom_\ZZ(P,\ZZ)$. We put
$\alpha_i := \eps_i - \eps_{i+1}$ ($1\leq i\leq n-1$).

The {\sl quantized universal enveloping algebra} $\Uq = 
\FSU_q(\gog\gol(n,\CC))$ (cf.\ \cite{\Dr}, \cite{\Jima}, \cite{\Jimb},
\cite{\Nmac})
 is the algebra generated by the symbols 
$q^h$ ($h\in P^\ast$) and $e_i, f_i$ ($1\leq i \leq n-1$) 
subject to the following relations:
$$
\eqalignno{ & q^0 =1, \quad q^{h+h'} = q^h\cdot q^{h'},&  \cr
& q^h e_i q^{-h} = q^{\langle h, \alpha_i\rangle} e_i, \quad
q^h f_i q^{-h} = q^{-\langle h, \alpha_i \rangle} f_i, & \cr
& e_if_j - f_j e_i = \delta_{ij} {q^{\alpha_i} -q^{-\alpha_i} 
\over q-q^{-1}}, & \eq{uqrel} \cr
& e_i^2 e_j - (q + q^{-1})e_ie_je_i + e_j e_i^2 = 0 \; 
(|i-j| = 1); \quad e_ie_j = e_je_i \, (|i-j| > 1), & \cr
&  f_i^2 f_j - (q + q^{-1})f_if_jf_i + f_j f_i^2 = 0 \; 
(|i-j| = 1); \quad f_if_j = f_jf_i \; (|i-j| > 1), & \cr}
$$
where $h,h'\in P^\ast$ and $1\leq i,j \leq n-1$.

Let $\Uh\subset \Uq$ denote the subalgebra generated by
the elements $q^h$ ($h\in P^\ast$). It is a Laurent polynomial
algebra in the generators $q^{\eps_i}$. 

A left $\Uq$-module $W$ is called {\sl $P$-weighted}
if it has a  vector space basis consisting of 
weight vectors with weights in $P$. The action of
$\Uh$ on any $P$-weighted $\Uq$-module is diagonalizable,
and any such module is completely reducible
(cf.\ \cite{\Ro}, \cite{\Lz}).
The cone of dominant weights $P^+\subset P$ consists of
all weights $\lambda = \sum_k \lambda_k\eps_k \in P$ such that 
$\lambda_1 \geq \ldots \geq \lambda_n$.
There is a 1-1 correspondence $\lambda \longleftrightarrow
V(\lambda)$ between dominant weights and irreducible
$P$-weighted finite-dimensional $\Uq$-modules 
such that $\lambda\in P^+$
is the highest weight of $V(\lambda)$
(cf.\ \cite{\Ro}, \cite{\Lz}).
Recall that 
$\lambda\in P^+$ is called a highest weight of an
irreducible  $\Uq$-module $V$
if there exists a non-zero vector $v\in V$ such that
$q^h\cdot v = q^{\langle h, \lambda\rangle} v$ and $e_i\cdot v = 0$
for all $1\leq i\leq n-1$.

We denote by $L^+_{ij}$, $L^-_{ij}$ the so-called {\sl L-operators}
introduced in \cite{\FRT}. For a precise definition of these operators
in the present context see \cite{\Nmac}. One has
$L^+_{ij} = L^-_{ji}=0$ ($i>j$). Moreover,
$$
q^h L^\pm_{ij} q^{-h} = q^{\langle h, \eps_j -\eps_i\rangle} 
L^\pm_{ij}, \quad
q^h S(L^\pm_{ij}) q^{-h} = q^{\langle h, \eps_j -\eps_i\rangle}
S(L^\pm_{ij}),
\eqno\eq{Lroot}
$$
for $1\leq i,j\leq n$. Because of this property,
the $L^\pm_{ij}$ can be viewed as $q$-analogues of the
root vectors in $\gog_\CC=\gog\gol(n,\CC)$.

The $L^\pm_{ij}$ generate the
algebra $\Uq$. There is a unique Hopf $\ast$-algebra structure
on $\Uq$ such that
$$
\Delta(L^\pm_{ij}) = \sum_k L^\pm_{ik} \ten L^\pm_{kj}, \quad
\eps(L^\pm_{ij}) = \delta_{ij}, \quad 
(L^\pm_{ij})^\ast = S(L^\mp_{ji}),
\eqno\eq{LDelta}
$$
for $1\leq i,j\leq n$. The subalgebra
$\Uh\subset \Uq$ is a Hopf $\ast$-subalgebra.
All $P$-weighted finite-dimensional representations of $\Uq$
are unitarizable.

Let $\rho_V\colon
\Uq \to \End(V)$ denote the irreducible representation
(called {\sl vector representation})
of $\Uq$ with highest weight $\eps_1\in P^+$. Put
$N:= \dim(V)$ (in the case $G=U(n)$, one actually has $N=n$).
Fixing  a suitable basis $(v_i)$ of $V$ such that $\Uh$ acts 
by diagonal matrices, one has
$$
R^\pm = \sum_{ij} e_{ij} \ten \rho_V(L^\pm_{ij})\in 
\End(V\ten V),
$$
where $R^\pm\in\End(V\ten V)$ 
(cf.\ \cite{\Jimc}, \cite{\FRT}, \cite{\Nmac})
 are the invertible $N^2\times N^2$ matrices defined by
$$
R := \sum_{ij} q^{\delta_{ij}} e_{ii}\ten e_{jj} +
(q-q^{-1}) \sum_{i>j} e_{ij}\ten e_{ji},\quad
R^+ := PRP, \quad R^- := R^{-1}.
\eqno\eq{Rdef}
$$
Here the $e_{ij}\in\End(V)$ are the standard 
unit matrices with respect to the basis $(v_i)$, and 
$P\in\End(V\ten V)$ is the usual permutation operator.

Denote by $\detq^{-1}\colon \Uq\to\CC$ the one-dimensional
representation with highest weight 
$\lambda = -\eps_1 - \cdots -\eps_n$.
The algebra $\Aq = \FSA_q(U(n))$  (cf.\ \cite{\FRT}, \cite{\Nmac}) of 
representative functions 
on the  {\sl quantum unitary group} $U_q(n)$ 
is the subalgebra of the linear dual of $\Uq$ generated by 
$\detq^{-1}$ and the coefficients $(t_{ij})$ 
of the vector representation $V$ with respect to the basis
$(v_i)$. The generators $t_{ij}$ and $\detq^{-1}$ satisfy the relations
$$
R T_1 T_2 = T_2 T_1 R, \qquad 
\detq \cdot \detq^{-1} = \detq^{-1} \cdot \detq.
\eqno\eq{RTT}
$$
Here $T := (t_{ij})_{1\leq i,j \leq n}$ is an 
$N\times N$ matrix with coefficients in $\Aq$, 
$T_1 := T\ten \id$ and 
$T_2 := \id\ten T$ are Kronecker matrix products, 
and $\detq\in \Aq$ 
(called the {\sl quantum determinant}) is defined as 
$$
\detq := \sum_{w\in\goS_n} (-q)^{l(w)} t_{w(1)1} \cdots 
t_{w(n)n}.
$$
Here $\goS_n$ denotes the permutation group on $n$ letters, 
and $l(w)$  denotes the length of a permutation 
$w\in \goS_n$. The elements $\detq$ and $\detq^{-1}$ are
 central in $\Aq$. 

There is a unique Hopf $\ast$-algebra structure on
$\Aq$ such that
$$
\Delta(t_{ij}) = \sum_k t_{ik} \ten t_{kj}, \quad \eps(t_{ij}) = 
\delta_{ij},\quad t_{ij}^\ast = S(t_{ji}),
$$
for all $1\leq i,j \leq n$, and
$$
\Delta(\detq) = \detq \ten \detq, \quad
\eps(\detq) =1, \quad (\detq)^\ast = S(\detq) = \detq^{-1}.
$$

The Hopf $\ast$-algebra $\Aq$
is by definition spanned by the coefficients of its 
finite-dimensional unitary
co\-representations (cf.\ \cite{\Wo}, \cite{\DKb}).
Any $P$-weighted representation of $\Uq$
can be lifted to a corepresentation
of $\Aq$.

Using the duality $\langle \cdot \, , \, \cdot \rangle$ 
between $\Uq$ and $\Aq$, 
one defines
a $\Uq$-bimodule structure on $\Aq$ by putting:
$$
u\cdot a := (\id \ten u) \circ \Delta(a), 
\quad a\cdot u := (u \ten \id) \circ \Delta(a),\quad 
u\in\Uq,\, a\in\Aq.
\eqno\eq{uqbimod}
$$
Here, on the right-hand side of both equalities,
$u\in \Uq$ is viewed as a linear form on $\Aq$ 
via the pairing $\langle \cdot \, , \, \cdot \rangle$. 
The multiplication
$\Aq \ten \Aq \to \Aq$ and the unit mapping $\CC \to \Aq$
then become $\Uq$-bimodule homomorphisms. In other words,
$\Aq$ is an {\sl algebra with two-sided $\Uq$-symmetry}.

Let $W(\lambda)\subset\Aq$ 
($\lambda\in P^+$) denote the subspace spanned by the coefficients of
the (co-)representation $V(\lambda)$. Then  one has the following 
multiplicity-free
decomposition of the $\Uq$-bimodule $\Aq$ into irreducible 
constituents:
$$
\Aq = \bigoplus_{\lambda\in P^+} W(\lambda), 
\quad W(\lambda) \simeq V(\lambda)^\ast \ten V(\lambda).
\eqno\eq{peterweyl}
$$
The above decomposition can also  be characterized as the simultaneous
eigenspace decomposition under the (left) action on $\Aq$ of the center 
$\FSZ\Uq\subset \Uq$.

Let $h\colon \Aq \to \CC$ denote the normalized {\sl Haar functional} 
(cf.\ \cite{\Wo}, \cite{\DKb}) on the compact quantum group
$G_q$. By putting $\langle a, b\rangle = h(b^\ast a)$ we define
a positive definite Hermitian form
$\langle \cdot \, , \, \cdot\rangle$ on $\Aq$. The subspaces
$W(\lambda)$ are mutually orthogonal with respect to this inner
product.

There is one particularly important subgroup inside the quantum
group $G_q$. 
Let $\AT:=\CC[z_1, \ldots, z_n, z_1^{-1}, \ldots, z^{-1}_n]$
be the algebra of Laurent polynomials in the variables $z_i$
($1\leq i\leq n$). 
We use the notation $z^\lambda := 
z_1^{\langle\lambda, \eps_1\rangle} 
\cdots z_n^{\langle\lambda, \eps_n\rangle}$ ($\lambda\in P$).
There is a unique Hopf $\ast$-algebra structure
on $\AT$ such that
$$ 
\Delta(z_i) = z_i \ten z_i, \quad \eps(z_i) = 1, \quad
z_i^\ast = z_i^{-1} \quad (1\leq i\leq n).
$$
The subalgebra $\Uh$ naturally is in Hopf $\ast$-algebra
duality with $\AT$ via
$$
\langle q^h, z^\lambda\rangle  := q^{\langle h, \lambda\rangle}.
$$
There is a unique surjective Hopf $\ast$-algebra morphism
$$
{}_{|\TT}\colon \Aq \longrightarrow \AT
\eqno\eq{toralrestr}
$$
mapping $t_{ij}\in \Aq$  onto $\delta_{ij} z_i\in \AT$ and 
$\detq\in \Aq$ onto $z_1 \cdots z_n\in \AT$. The torus
$\TT$ can be viewed as a maximal toral subgroup of the quantum 
group $G_q$. The mapping ${}_{|\TT}$ then 
is the corresponding restriction of functions.
The mappings ${}_{|\TT}$ and $\Uh \hookrightarrow \Uq$ are dual
to each other. 

In the classical case, one can freely move around the maximal
torus $\TT\subset G$ by conjugation. After quantization,
however, this is no longer true. Since 
the group of Hopf algebra automorphisms of $\Aq(G)$ 
 is rather small (cf.\ \cite{\Chin}),
there is very little possibility for changing the
position of $\TT$ inside the quantum group $G_q$. This
fact will prove to be of cardinal importance in the
analysis of quantum symmetric spaces.
\beginsektion 3 {Quantum homogeneous spaces}
We recall some general facts (cf.\ \cite{\DKa}).
The ``infinitesimal'' method to construct quantum homogeneous spaces
uses the idea that the algebra of functions
on the homogeneous space $G/K$ can also be defined as the subspace
of those functions on $G$ which are annihilated by the $G$-invariant
differential operators $X\in \gok_\CC\subset \gog_\CC$. 

In fact,
suppose that $\gok_q \subset \Uq$ is a two-sided coideal invariant
under the mapping $\tau=\ast\circ S\colon \Uq\to \Uq$. 
Recall that a subspace $\gok_q\subset  \Uq$ is called a two-sided coideal if
$\Delta(\gok_q)\subset \gok_q\ten\Uq + \Uq\ten\gok_q$ and $\eps(\gok_q)=0$.
Then the subspace
$$
B_{\gok_q} := \{a\in\Aq \mid a\cdot \gok_q = 0\}
\eqno\eq{Bkdef}
$$
of right $\gok_q$-invariant functions
is a $\ast$-subalgebra and right coideal in $\Aq$, 
and a left $\Uq$-submodule.

Suppose we have a surjective Hopf $\ast$-algebra morphism
$\pi\colon \Aq(G) \to \Aq(K)$ and a dual Hopf $\ast$-algebra
mapping $\psi\colon \Uq(\gok_\CC)\to\Uq(\gog_\CC)$. Then
the subspace $\gok_q\subset\Uq(\gog_\CC)$ spanned by the elements
$\psi(v)-\eps(v) 1$ ($v\in \Uq(\gok_\CC)$) is a $\tau$-invariant
coideal. Left or right invariance with respect to the
quantum subgroup $K_q\hookrightarrow G_q$ then is the same as
invariance with respect to $\gok_q$. 
Unfortunately, as is well-known,
the quantization procedure as given by Drinfeld and Jimbo
is not functorial (cf.\ \cite{\Hay}). 
In other words, to an embedding of compact 
Lie groups $K\hookrightarrow G$ there need not correspond a
surjective Hopf $\ast$-algebra morphism $\Aq(G) \to \Aq(K)$.
Actually, 
the supply of Riemannian symmetric pairs $(G,K)$ whose embedding
$K\hookrightarrow G$ survives quantization turns out to be rather
limited (cf.\ section 5).

Suppose we have a $\tau$-invariant two-sided coideal 
$\gok_q\subset \Uq$. Then the pair $(\Uq, \gok_q)$ is called a
{\sl quantum Gelfand pair} if, for any finite-dimensional
irreducible $P$-weighted representation $V$ of $\Uq$,
the subspace $V_{\gok_q} := \{v\in V\mid \gok_q \cdot v=0\}$ of
so-called {\sl $\gok_q$-fixed vectors} is at most one-di\-men\-sion\-al.
Those representations $V$ for which the subspace $V_{\gok_q}$ is 
non-zero are called {\sl spherical}. Let us denote the
corresponding subset  of dominant weights by $P^+_{\gok_q}$.
Given a quantum Gelfand pair $(\Uq, \gok_q)$ one can define the
$\ast$-algebra $\FSH$ of $\gok_q$-biinvariant functions as
$$
\FSH := \{a\in \Aq\mid a\cdot \gok_q = 0 \; \hbox{and}\; \gok_q\cdot
a = 0 \}.
\eqno\eq{Hdecomp}
$$
Once again, one has a canonical decomposition
$\FSH = \bigoplus_{\lambda\in P^+_{\gok_q}} \FSH(\lambda)$,
where $\FSH(\lambda)$ is the intersection of $\FSH$ with
the subspace $W(\lambda)\subset \Aq$. 
Each of the subspaces
$\FSH(\lambda)$ is one-di\-men\-sion\-al. One can now define a 
{\sl (zonal) spherical function} as a non-zero element
of $\FSH(\lambda)$ ($\lambda\in P^+_{\gok_q}$). Spherical functions
corresponding to different $\lambda\in P^+_{\gok_q}$ will be orthogonal
with respect to the inner product on $\Aq$ defined in terms
of the Haar functional $h\colon \Aq \to \CC$.

Suppose now 
that the restriction of the mapping
${}_{|\TT}\colon \Aq \to \AT$ to $\FSH$ is injective onto its
image $\FSH_{|\TT}$ (as we will see, there are interesting cases
in which this condition is not fulfilled). 
The algebra $\FSH$ will  be
commutative then. Let 
$C\in\FSZ\Uq$ be a suitable central element ({\sl Casimir operator}).
The left action of $C$ on $\Aq$ will preserve the subalgebra
$\FSH$. Hence, there is a uniquely determined operator 
$D\colon \FSH_{|\TT}\to \FSH_{|\TT}$ such that on
$\FSH$ we have 
$$
{}_{|\TT} \circ C = D \circ {}_{|\TT},
\eqno\eq{raddef}
$$
where the symbol $C$ denotes the left action of
$C\in \Uq$ on the subalgebra $\FSH$. The operator $D$ will
be called the {\sl radial part} of the Casimir operator $C$.
Since $C$ is central, it acts as a scalar on
every subspace $\FSH(\lambda)$ ($\lambda\in P^+_\gok$).
In other words, the restriction $\phi(\lambda)_{|\TT}$ 
of  the spherical function $\phi(\lambda)\in\FSH(\lambda)$ 
($\lambda\in 
P^+_\gok$) to the maximal torus $\TT$ is an eigenfunction
of the operator $D$. If we are able to compute  an explicit
expression for the operator $D$, this is likely to give us
a strong clue to the nature of the spherical function
$\phi(\lambda)$.

Summarizing, we can say that one possible 
approach to tackling the problem
of quantizing the symmetric pair $(G,K)$ and describing 
the corresponding $q$-spherical functions consists of
the following steps. 
{\parindent=40pt
\item{(i)}
Find an
appropriate $\tau$-invariant coideal $\gok_q\subset \Uq$ which is a 
$q$-analogue of the Lie algebra $\gok_\CC$ of $K$. Show
that the pair $(\Uq, \gok_q)$ is a quantum Gelfand pair and 
that the spherical representations are indexed by the same
subset of highest weights as in the classical case.
\item{(ii)}
Investigate whether $\gok_q$-biinvariant
functions on $G_q$ are completely determined by their restriction to
the maximal torus $\TT\subset G_q$, and, if so, give an explicit
description of the image $\FSH_{|\TT}$ of $\FSH$ under
the restriction mapping ${}_{|\TT}$. 
\item{(iii)}
Find an explicit expression for the radial part $D$ of
some suitable Casimir element $C\in\FSZ\Uq$.
\par}
\beginsektion 4 {Quantum projective spaces} 
In this section we carry out the programme sketched at the end of section
3 for the case of a complex projective space.
The results in this section are joint work with M.\ Noumi. For proofs
and more details the reader is referred to the forthcoming paper \cite{\DN}.

We briefly review the classical case to provide motivation for what follows.
Set $G=U(n)$, $K=U(n-1)\times U(1)\subset U(n)$. The corresponding
involution $\theta\colon \gog_\CC \to\gog_\CC$ is given by
$\theta(X) = J XJ^{-1}$ with
$$
J= \hbox{diag}\{1, \ldots, 1, -1\}.
\eqno\eq{Jdef}
$$
One has $P_K^+ = \{l(\eps_1 -\eps_n)\in P\mid l\in \ZZ_+\}$.
The Lie subalgebra $\gok_\CC$ is equal to $\gog\gol(n-1,\CC)\oplus
\gog\gol(1,\CC)$. The maximal abelian
subspace $\goa_\CC\subset \gop_\CC$ is the one-dimensional subspace
spanned by the matrix $X=e_{1n} + e_{n1}$. The restricted root system
is isomorphic with $BC_1$. The root multiplicities are $1$ (long roots)
and $2(n-2)$ (short roots). 

Let us recall that the (classical) {\sl Jacobi poynomials} 
$P^{(\alpha, \beta)}_n(x)$ are one-variable  polynomials which are
orthogonal on the interval $[-1,1]$ 
with respect to the continuous weight function 
$w(x) = (1-x)^\alpha(1+x)^\beta$ ($\alpha, \beta > -1$). 
They actually coincide with the generalized Jacobi polynomials
corresponding to the root system $BC_1$ (there is a simple linear
relation between the parameters $\alpha$, $\beta$, and the values of
the multiplicity function $k$).
The usual Casimir operator $C\in\FSU(\gog_\CC)$ induces a 
$G$-invariant differential operator on $G/K$ whose radial part 
essentially coincides with the differential operator 
in \eqtag{jacdiff} for the values $\alpha = n-2$,
$\beta=0$. Hence, the spherical function $\phi_l$ corresponding
to the highest weight $l(\eps_1-\eps_n)$ $(l\in\ZZ_+$) can be expressed
as a Jacobi polynomial $P^{(n-2,0)}_l$.

Our choice of the involution $\theta$  is such that the corresponding
maximal abelian subspace $\goa_\CC$ does not contain any diagonal
matrices. Hence the restriction of $K$-biinvariant functions to the
diagonal subgroup $\TT$ is hardly injective.
To remedy this situation, one could take a different choice of
involution, for instance by putting
$$
J' := \hbox{diag}(0,1,\ldots, 1,0) -e_{n1} -e_{n1},
\eqno\eq{Jprimedef}
$$
and defining $\theta'(X) = J'X{J'}^{-1}$. The corresponding
$\goa'_\CC$ is spanned by the matrix $X'=e_{11}-e_{nn}$ and
hence is contained in the Cartan subalgebra $\goh_\CC$.
With this choice of $\theta'$ and the corresponding subgroup 
$K'\subset G$, the restriction mapping $\FSH'=\FSA(K'\backslash G/K') \to
\FSA(\TT)$ actually is injective. Of course, the symmetric spaces
$G/K$ and $G/K'$ are isometric and their spherical functions
essentially the same, since the involutions $\theta$ and $\theta'$
differ only by a conjugation. Alternatively, one could
have moved the maximal torus $\TT$ inside $G$ so as to make
the restriction of $K$-biinvariant functions injective.

We now turn to the quantum case.
Let us fix real numbers $c,d\geq 0$ such that 
$(c,d) \neq (0,0)$. 
The subspace 
$\gok^{(c,d)} \subset \Uq$ is by 
definition spanned 
by the following elements:
$$
\eqalignno{(i) \; & \Lp_{11} - \Lm_{nn}, \;
\Lm_{11} - \Lp_{nn}, & \cr
(ii) \; & \sqrt{c}\, \Lp_{1k} + \sqrt{d}\, 
\Lm_{nk} \quad (2\leq k\leq n-1), & \cr
(iii) \; & \sqrt{d}\, \Lp_{kn} + 
\sqrt{c}\, \Lm_{k1} \quad
(2\leq k\leq n-1),& \eq{kdef} \cr
(iv) \; & \Lp_{ij}, \; \Lm_{ji} \quad 
(2\leq i<j\leq n-1),& \cr
(v) \; &\Lp_{ii} - \Lm_{ii} \quad 
(2\leq i \leq n-1), & \cr
(vi) \; & \sqrt{cd}\, \Lp_{1n} - 
\sqrt{cd}\, \Lm_{n1}
-(c-d) (\Lp_{11} -\Lm_{11}). & \cr}
$$
We remark that the subspace $\gok^{(c,d)}\subset \Uq$ only 
depends on the ratio of the numbers $c$ and $d$. In fact,
given $c,d\geq 0$ such that 
$(c,d) \neq (0,0)$, we can define $\sigma\in\RR \cup \{\pm \infty\}$ by
$q^{\sigma} = \sqrt{{d\over c}}$ ($c,d>0$), 
$\sigma=-\infty$ ($c=0$), $\sigma=\infty$ ($d=0$).
Then the coideal $\gok^{(c,d)}$ only depends
on the value of the corresponding $\sigma$. We write 
$\gok^\sigma = \gok^{(c,d)}$. 

In case $\sigma$ is finite, there is the following alternative, 
but equivalent, way to 
define the subspace $\gok^\sigma$. Let $J^\sigma\in \End(V)$ 
be the $n\times n$ matrix with 
complex coefficients defined by:
$$
J^\sigma := \hbox{diag}(q^\sigma(q^{-\sigma} -q^\sigma), 1,
\ldots, 1,0) -q^\sigma e_{1n}-q^\sigma e_{n1}.
\eqno\eq{Jsigmadef}
$$
Let $M^\sigma \in \End(V)\ten\Uq$ be the $n\times n$ matrix with
coefficients in $\Uq$ defined by 
$M^\sigma := L^+ J^\sigma -J^\sigma L^-$. Then $\gok^\sigma$ is
spanned by the coefficients $M^\sigma_{ij}$ of the matrix $M^\sigma$.

In the limit $q\to 1$, the subspaces $\gok^\infty$ and
$\gok^{-\infty}$ will essentially
tend to the Lie subalgebras $\gog\gol(n-1,\CC)\oplus \gog\gol(1,\CC)$
and $\gog\gol(1,\CC)\oplus \gog\gol(n-1,\CC)$ respectively.
Moreover, the subspace $\gok^0$ will essentially tend to the Lie subalgebra
$\gok'_\CC$ defined with respect to \eqtag{Jprimedef}.
\beginth{Lemma} kcoid
{\thstil Suppose $-\infty \leq \sigma \leq \infty$.
The subspace $\gok^\sigma\subset \Uq$ is a 
$\tau$-invariant two-sided coideal in $\Uq$.
\par}
\endth
There are natural surjective Hopf $\ast$-algebra morphisms
$$
\Aq \to \Aq(U(n-1))\ten \FSA(U(1)), \quad
\Aq\to \FSA(U(1))\ten \Aq(U(n-1))
$$
corresponding to the embeddings $U(n-1)\times U(1) \subset U(n)$ etc.
Left or right invariance with respect to these quantum subgroups
is the same as invariance with respect to the
coideals $\gok^\infty$ and $\gok^{-\infty}$ respectively.
For finite values of $\sigma$, it can be shown that 
there is no quantum subgroup
corresponding to the coideal $\gok^\sigma$.
\beginth{Lemma} hweightcomp
{\thstil Suppose $-\infty  <\sigma < \infty$. 
Let $V$ be a finite-dimensional 
irreducible $P$-weighted $\Uq$-module. If $v\in V$ is 
a non-zero $\gok^\sigma$-fixed vector then the highest weight 
component of $v$ is non-zero.
\par}
\endth
\beginth{Theorem} spherrep
{\thstil Suppose $-\infty \leq \sigma \leq \infty$. 
For all $\lambda \in P^+$, the subspace $V(\lambda)_{\gok^\sigma}$ of
$\gok^\sigma$-fixed vectors is at 
most one-di\-men\-sion\-al. The subspace
$V(\lambda)_{\gok^\sigma}$ is one-di\-men\-sion\-al if and only if 
$\lambda = l(\eps_1 - \eps_n)$ for some $l\in \ZZ_+$.
\par}
\endth
In other words, the spherical representations are labelled by the 
same highest weights as in the classical case.
We put
$$
\Bq^{(c,d)}:= \{ a \in \Aq \mid a \cdot \gok^{(c,d)} = 0 \}.
$$
 The subspace
$\Bq^{(c,d)}$ is a $\ast$-subalgebra and right coideal in $\Aq$,
and it is invariant under the left action of $\Uq$ on $\Aq$.
Vaksman and Korogodsky \cite{\VK} defined the
 $\Aq$-comodule algebra $\Bq^{(c,d)}$ (arbitrary values of $c,d$)
by means of a $q$-analogue of the Hopf fibration
$S^{2n-1}\to \CC\PP^{n-1}$. It can be shown that the comodule algebras
$\Bq^\sigma = \Bq^{(c,d)}$ are non-isomorphic for different values
of $\sigma$.
\beginth{Theorem} planch
{\thstil Suppose $c,d\geq 0$ such that $(c,d) \neq (0,0)$.
The irreducible decomposition of $\Bq^{(c,d)}$ as a
right $\Aq$-comodule resp.\ left $\Uq$-module is given by:
$$
\Bq^{(c,d)} = \bigoplus_{l\in\ZZ_+} V(l(\eps_1 - \eps_n)),
$$
where the isotypical subspace of type $V(l(\eps_1 - \eps_n))$ 
is equal to the intersection of $\Bq^{(c,d)}$ and $W(\lambda)$.
\par}
\endth
We put
$$
z_{ij} := d t_{1i}^\ast t_{1j} + c t_{ni}^\ast t_{nj}
+ \sqrt{cd}\, t_{ni}^\ast t_{1j} + \sqrt{cd}\, t_{1i}^\ast t_{nj}
\in\Aq \quad (1\leq i,j \leq n).
\eqno\eq{zdef}
$$
\beginth{Proposition} zgenB
{\thstil Suppose $c,d\geq 0$ such that $(c,d)\neq (0,0)$.
The $z_{ij}$ ($1\leq i,j \leq n$) are right $\gok^{(c,d)}$-invariant
 and generate the subalgebra $\Bq^{(c,d)}$. They satisfy
$z_{ij}^\ast = z_{ji}$.
\par}
\endth
Let $V^\ast$ denote the contragredient of the 
vector representation $V$. It has highest weight $-\eps_n$.
Let $(v^\ast_i)$ denote the dual basis of $V^\ast$. The tensor
product representation $V^\ast \ten V$ has the irreducible
decomposition 
$$
V^\ast \ten V \cong V(0) \oplus V(\eps_1 - \eps_n).
$$ 
The subspace of $\gok^{(c,d)}$-fixed vectors in $V^\ast \ten V$ is
two-dimensional and spanned by the two elements
$$
\sum_k q^{2(n-k)} v^\ast_k\ten v_k, \quad
\sqrt{cd}\,v^\ast_1 \ten v_n + \sqrt{cd}\, v^\ast_n \ten v_1 + 
 q d\, v^\ast_1 \ten v_1 + q^{-1} c\, v^\ast_n \ten v_n.
\eqno\eq{kfixedV}
$$

Let us fix parameters $-\infty\leq\sigma,\tau\leq \infty$. 
We put
$$
\Hst := \{ a\in \Aq \mid
\gok^\sigma \cdot a = 0 \; \hbox{and} 
\; a\cdot \gok^\tau = 0\},
\eqno\eq{Hstdef}
$$
and call it the space of {\sl $(\sigma,\tau)$-biinvariant}
functions. It actually is a $\ast$-subalgebra of $\Aq$.
If we put $\Hst(\lambda) := \Hst\cap W(\lambda)$, then 
$$
\Hst = \bigoplus_{l\in\ZZ_+} \Hst(l(\eps_1 -\eps_n)).
$$
Each of the spaces $\Hst(l(\eps_1 -\eps_n))$ is
one-dimensional. A {\sl $(\sigma,\tau)$-spherical function}
is by definition a non-zero element of $\Hst(l(\eps_1-\eps_n))$
($l\in\ZZ_+$).

Let us now suppose that $\sigma,\tau$ are finite.
The direct sum $\Hst(0) \oplus \Hst(\eps_1
-\eps_n)$ is spanned by the unit element $1\in \Aq$
and the element (cf.\ \eqtag{kfixedV})
$$
\xst := {1\over 2} (z_{1n} + z_{n1} + q^{\sigma +1}z_{11} 
 + q^{-\sigma-1} z_{nn} -
(q^{\sigma+\tau+1} + q^{-\sigma-\tau-1}))
\in \Aq. 
\eqno\eq{xstdef}
$$
Here (cf.\ \eqtag{zdef})
$$
z_{ij} := q^\tau t_{1i}^\ast t_{1j} + 
q^{-\tau} t_{ni}^\ast t_{nj}
+ t_{ni}^\ast t_{1j} + t_{1i}^\ast t_{nj}
\in\Aq \quad (1\leq i,j \leq n). 
$$
Note that $(\xst)^\ast = \xst$.
\beginth{Lemma} Trestrlem
{\thstil Let $\sigma,\tau$ be finite. 
Under the restriction mapping 
${}_{|\TT}\colon \Aq \longrightarrow \AT$ we have:
$$
z_{11} \mapsto q^\tau, \quad z_{nn} 
\mapsto q^{-\tau}, \quad
z_{1n} \mapsto z_1^{-1} z_n, \quad z_{n1} 
\mapsto z_1z_n^{-1},
$$
and all the other $z_{ij}$ are mapped onto $0$. 
In particular, the image of $\xst$ equals 
${1\over 2}(z_1z_n^{-1} + z_1^{-1}z_n)$.
\par}
\endth
\beginth{Theorem} Trestrinj
{\thstil Let $\sigma,\tau$ be finite.
We put $z := z_1z_n^{-1}\in \AT$.
The restriction of the mapping ${}_{|\TT}$ to 
$\Hst\subset \Aq$
is an injective $\ast$-algebra homomorphism 
onto the polynomial algebra 
$\Hst_{|\TT}:=\CC[{1\over 2} (z+z^{-1})] 
\subset \AT$.  Hence, the algebra
$\Hst= \CC[\xst]$ is commutative.
\par}
\endth

For every $l\in\ZZ_+$, let us fix a non-zero 
$\phi_l\in\Hst(l(\eps_1-\eps_n))$. The 
$(\sigma,\tau)$-spherical function  $\phi_l$ is
uniquely determined up to multiplication by a scalar.
 It can be expressed as a polynomial of 
degree $l$ in $\xst$.  In order to identify 
these polynomials, we study the action of
the following {\sl Casimir operator} on $\Hst$
(cf.\ \cite{\FRT}, \cite{\Nmac}):
$$
C := \sum_{ij} q^{2(n-i)}\Lp_{ij}S(\Lm_{ij})\in 
\Uq.
\eqno\eq{Cdef}
$$
The element $C$ acts as a scalar on each subspace 
$\Hst(l(\eps_1-\eps_n))$ $(l\in \ZZ_+$). 
The corresponding eigenvalue is
$$
q^{2(l+n-1)} + q^{-2l} + \sum_{1<i<n} q^{2(n-i)}.
$$
As in \eqtag{raddef}, we can consider the radial part
$D\colon \Hst_{|\TT} \to \Hst_{|\TT}$ of the Casimir
operator $C$.
Let us define a linear operator $T_{q,z} \colon
\CC[z,z^{-1}] \to \CC[z,z^{-1}]$ by 
putting $T_{q,z} f(z) := f(qz)$. 
We use the following notation:
$$
A(z;q) := 
{(1-az)(1-bz)(1-cz)(1-dz)\over (1-z^2)(1-qz^2)} \quad
(a,b,c,d\in\CC).
\eqno\eq{Azqdef}
$$
\beginth{Theorem} radialpart
{\thstil Let $\sigma,\tau$ be finite. 
The radial part $D\colon \Hst_{|\TT} \to \Hst_{|\TT}$ 
of the Casimir operator $C$ is
equal to the following second-order $q$-difference
operator:
$$
D= A(z;q^2)(T_{q^2,z} - \id)  + 
A(z^{-1};q^2)(T_{q^{-2},z} -\id)
+ {1-q^{2n} \over 1-q^2}\cdot\id,
$$
with parameters $a,b,c,d$ given by
$$
a = -q^{\sigma + \tau + 1},
\quad b= -q^{-\sigma - \tau + 1},
\quad c = q^{\sigma - \tau + 1}, \quad
d = q^{-\sigma + \tau + 2(n-2) + 1}.
$$
\par}
\endth
\noindent
For $n=2$ this result was essentially proved in Koornwinder
\cite{\KWd, Lemma 5.2}.

Let us recall the definition of Askey-Wilson 
polynomials (cf.\ \cite{\AW}).
The  {\sl $q$-shifted factorials}
are defined as
$$
(a;q)_n := \prod_{k=0}^{n-1} (1-aq^k), \;
(a_1,\ldots, a_s;q)_n := \prod_{j=1}^s (a_j;q)_n,
\; (a;q)_\infty := \lim_{n\to\infty} (a;q)_n,
$$
and the {\sl $q$-hypergeometric series} (cf.\ \cite{\GR}) as
$$
{}_{s+1}\phi_s \left [ {a_1, \ldots, a_{s+1} \atop
b_1, \ldots, b_s}; q,z \right ] :=
\sum_{k=0}^\infty {(a_1, \dots, a_{s+1};q)_k \,z^k\over
(b_1, \ldots, b_s;q)_k\, (q;q)_k}.
$$
{\sl Askey-Wilson polynomials}
are defined as the polynomials in $\cos(\theta)$
given by
$$
\eqalignno{p_n(\cos(\theta)&; a,b,c,d \mid q) := & \cr
& a^{-n} (ab,ac,ad;q)_n  \cdot
{}_4\phi_3\left [ {q^{-n}, q^{n-1}abcd, ae^{i\theta}, 
ae^{-i\theta}\atop
ab,ac,ad}; q,q \right ]. & \eq{awdef}\cr}
$$
They are symmetric in the parameters $a,b,c,d$. 
Suppose that $a,b,c,d$ are real, or if complex, appear in
complex conjugate pairs, and that $|a|, |b|, |c|, |d|
\leq 1$ such that  the pairwise products of $a,b,c,d$ are
not $\geq 1$. Then the Askey-Wilson polynomials
$p_n$ satisfy the orthogonality relations (cf.\ \cite{\AW}))
$$
\int_0^{2\pi} (p_np_m)(\cos \theta; a,b,c,d \mid q)
\left \vert {(e^{2i\theta};q)_\infty\over 
(ae^{i\theta}, be^{i\theta}, ce^{i\theta}, de^{i\theta};q)_\infty}
\right \vert^2 \hbox{d}\theta = 0, \; n\neq m.
\eqno\eq{aworth}
$$
If the condition 
$|a|, |b|, |c|, |d| \leq 1$ is relaxed, finitely many discrete
terms will appear in the orthogonality relation \eqtag{aworth}.

Askey-Wilson polynomials can also be 
characterized as the eigenfunctions of a certain 
second-order $q$-difference operator. To be more precise,
if we write the ${}_4\phi_3$-factor in \eqtag{awdef} as 
$r_n(\cos(\theta))$ and put 
$P_n(z) := r_n({1\over 2}(z+z^{-1}))$,
then the Laurent polynomials $P_n(z)$ in the variable $z$
satisfy the following second-order $q$-difference
equation:
$$
\eqalignno{A(z;q)(P_n(qz) - P_n(z))&  + 
A(z^{-1};q)(P_n(q^{-1}z) -P_n(z)) = & \cr
& = -(1-q^{-n})(1-q^{n-1}abcd)P_n(z), & \eq{awqdiff}\cr}
$$
where $A(z;q)$ is defined in \eqtag{Azqdef}.
Any polynomial $f(e^{i\theta})$ of degree $\leq n$ 
in $\cos \theta$ satisfying \eqtag{awqdiff}
with $z=e^{i\theta}$ is a constant
multiple of $P_n(e^{i\theta})$.

Comparing \eqtag{awqdiff} with \thtag{radialpart} we
conclude:
\beginth{Theorem} awspher
{\thstil Let $\sigma,\tau$ be finite. The 
$(\sigma,\tau)$-spherical functions in  
$\Hst(l(\eps_1-\eps_n))$ $(l\in\ZZ_+)$ 
are spanned by
$$
p_l(\xst; -q^{\sigma + \tau + 1}, 
-q^{-\sigma - \tau + 1},
q^{\sigma - \tau + 1}, 
q^{-\sigma + \tau + 2(n-2) + 1} \mid q^2),
$$
where $p_l$ is an Askey-Wilson polynomial.
\par}
\endth\noindent
For $n=2$ this result  was essentially proved by
Koornwinder \cite{\KWd, Theorem 5.2}.

Using \thtag{awspher} and a suitable 
limit argument (cf.\ \cite{\KWd, Prop.\ 6.1, 6.3}), 
one can identify the spherical
functions corresponding to the cases when either 
$\sigma$ or $\tau$ is infinite. They are
expressed as so-called {\sl little} and {\sl big 
q-Jacobi polynomials} (cf.\ \cite{\AA}). 
These are orthogonal polynomials
with respect to an measure supported on
an infinite discrete set. For more details see
\cite{\DN}.
\beginsektion 5 {Higher rank quantum symmetric spaces} 
We refer to Table I for an exhaustive list (up to local
isomorphism) of the classical 
irreducible compact Riemannian symmetric spaces of
type I (cf.\ \cite{\Hela}). The column labelled CC
contains the type designation as used 
by \'E.\ Cartan in his classification. In the last column
we give the multiplicities of the short, medium-sized, and
long restricted roots (in this order). 
Here we consider $C_l$ as a root system
of type $BC$, the short roots having multiplicity zero. 
\bigskip
\centerline{Table I: Classical irreducible compact symmetric spaces}
\medskip
\hrule
\settabs \+ xxx &  AIIII\ \ & $SO(2n+1)$\ \  & $SO(l)\times SO(2n+1-l)$
\ \ & $l=2n+1$\ \ & $BC_l$\ \ & xxxxxx \cr
\+ & & & & & & \cr
\+ no.\ & CC & $G$ & $K$ & $l$ & $\Sigma$ & $m_\alpha$ \cr
\+ & & & & & & \cr
\hrule
\+ & & & & & & \cr
\+ 1 & AI & $SU(n)$ & $SO(n)$ & $l=n-1$ & $A_l$ & 1\cr
\+ 2 & AII & $SU(2n)$ & $Sp(n)$ & $l=n-1$ & $A_l$ & 4\cr
\+ 3& AIII & $U(n)$ & $U(l)\times U(n-l)$ & $l\leq \left [{n\over 2}\right ]$
 & $BC_l$ & $2(n-2l), 2, 1$ \cr
\+  4 & BI & $SO(2n+1)$ & $SO(l)\times SO(2n+1-l)$ & $l\leq n$ & 
$B_l$ & $1, 2n+1-2l$\cr
\+ 5 & CI & $Sp(n)$ & $U(n)$ & $l=n$ & $C_l$ & $0,1,1$\cr
\+ 6 & CII & $Sp(n)$ & $Sp(l)\times Sp(n-l)$ & $l\leq \left [ {n\over 2}
\right ]$ & $BC_l$ & $4(n-2l), 4, 3$\cr
\+ 7 & DI & $SO(2n)$ & $SO(l)\times SO(2n-l)$ & $l < n$ & $B_l$ &
$2n-2l, 1$\cr
\+ 8 & DI & $SO(2n)$ & $SO(l)\times SO(2n-l)$ & $l = n$ & $D_l$ &
1\cr
\+ 9 & DIII & $SO(2n)$ & $U(n)$ & $n=2l$ & $C_l$ & $0,4,1$\cr
\+ 10 & DIII &  $SO(2n)$ & $U(n)$ & $n=2l+1$ & $BC_l$ & $4,4,1$\cr
\+ & & & & & & \cr
\hrule\bigskip
Let $(G,K)$ be any of the symmetric pairs listed in Table I.
Suppose we have a two-sided coideal $\gok_q\subset \Uq$ which,
in some suitable sense, has $\gok_\CC \subset \gog_\CC$ as
its classical limit.
From the results in section 4 we may draw the conclusion
that the relative position of $\gok_\CC$ with respect to the
standard Cartan subalgebra $\goh_\CC\subset \gog_\CC$ 
is of considerable importance in the analysis of the quantum
symmetric space corresponding to $\gok_q$ and its zonal
spherical functions. We distinguish three different cases.

The first case arises when $\goh_\CC$ is $\theta$-invariant and 
$\gok_\CC$ such that the intersection
$\goh_\CC \cap \gop_\CC$ is maximal abelian in $\gop_\CC$, hence
equal to $\goa_\CC$. We shall term this the {\sl maximally split}
case. In this situation one can reasonably expect the restriction
of $\gok_q$-biinvariant functions to the diagonal subgroup $\TT$
to be injective and the spherical functions to be multivariable
orthogonal polynomials with respect to a {\sl continuous} measure.

Actually, so far, Noumi and Sugitani \cite{\Nmac}, \cite{\NSa},
\cite{\NSb}, \cite{\Su} have constructed a maximally split coideal
$\gok_q\subset \Uq$ for all symmetric pairs
$(G,K)$ listed in Table I, and analysed the corresponding
quantum symmetric spaces  and their spherical functions.
We shall give a very brief account of their results
for the symmetric pairs (1), (2), and (4)--(10).

These symmetric pairs have in common that the ``lowest''
non-trivial spherical representation occurs in the irreducible
decomposition of the tensor product $V\ten V$ of the vector 
representation $V$ with itself. In the quantum case, the
existence of $\gok_q$-fixed vectors inside $V\ten V$ is controlled
by the so-called {\sl reflection equation}
(cf.\ \cite{\Ku}). Suppose we have an invertible complex $N\times N$ matrix
$J$ (recall that $N := \dim(V)$). Define a matrix
$M\in \End(V)\ten \Uq$ with coefficients in $\Uq$ by
putting 
$$
M:= L^+ - J S(L^-)^t J^{-1},
\eqno\eq{MJdef}
$$
and define $\gok_q\subset \Uq$ as the subspace spanned by the
coefficients $M_{ij}$ of $M$.  Then $\gok_q$ will be a two-sided
coideal in $\Uq$. Put $w_J = \sum_{ij} v_i\ten J_{ij}v_j \in V\ten V$.
Then it can be shown that $\gok_q \cdot w_J = 0$ if and only if
the matrix $J$ satisfies the following reflection equation:
$$
R J_1 R^{t_1} J_2 = J_2 R^{t_1} J_1 R,
\eqno\eq{reflVV}
$$
where ${}^{t_1}$ denotes transposition in the first tensor component, 
and $J_1 := J\ten \id$ etc.\ are Kronecker matrix products.

For every symmetric pair $(G,K)$, Noumi and Sugitani give 
a constant (i.e.\ involving no parameters besides $q$) matrix $J$ 
satisfying the reflection equation and such that
the corresponding $\gok_q$ is invariant under $\tau=\ast\circ S$.
Let $\rho\in P$ be the half sum of the positive roots. The element
$q^\rho\in\Uh$ acts on $V$ as a diagonal matrix written
$\hbox{diag}(q^{\rho_1}, \ldots, q^{\rho_N})$. 
Then the matrix $J$ is given by:
$$
\eqalign{J &=  \hbox{diag}(q^{\rho_1}, \ldots, q^{\rho_N})\quad
\hbox{for} \quad (1),\, (5),\, (8),\cr
J & =  J_0\hbox{diag}(q^{\rho_1}, \ldots, q^{\rho_N}) \quad
\hbox{for} \quad (2),\, (6)\, (n=2l), (9), (10),\cr}
$$
where $J_0 := \sum_{k=1}^N (-e_{2k,2k-1} + e_{2k-1,2k})$. For
the definition of $J$ in the remaining cases see \cite{\NSb}, 
\cite{\Su}.
In all cases,  the pairs $(\Uq, \gok_q)$ are quantum Gelfand pairs.
The spherical representations are labelled by the same subset
$P^+_K$ of highest weights as in the classical case.

Noumi and Sugitani next study the $\ast$-subalgebra $\FSH\subset \Aq$
of $\gok_q$-biinvariant functions. To be precise, in the cases
(4)--(10) they consider functions which are left-invariant
with respect to $\gok_q$ and right invariant with respect to
$\overline{\gok_q}$, where $\bar u := q^\rho u^\ast q^{-\rho}$ is
a modified $\ast$-operation on $\Uq$. 
They show that the mapping ${}_{|\TT}\colon \FSH \to \FSA(\TT)$ is
injective and  describe the image $\FSH_{|\TT}\subset \FSA(\TT)$
in terms of certain explicitly defined
elements $x_1, \ldots , x_l \in \FSA(\TT)$. More precisely, 
they show this subalgebra
is isomorphic with $\FSA(\Sigma)^W$ under the
correspondence $x_i\mapsto e^{\eps_i}$. 

Let us recall the definition of Macdonald's symmetric orthogonal
polynomials (cf.\ \cite{\Ma}, \cite{\Mb}, \cite{\Mc}) and Koornwinder's
$BC$-type Askey-Wilson polynomials (cf.\ \cite{\KWc}). In Macdonald's
case, let $\Sigma$ be one of the root systems $A_l$, $B_l$,
$C_l$, $D_l$. In 
Koornwinder's case we put $\Sigma=BC_l$. 
We freely use the notation introduced in section 1.

We first treat Macdonald's case. 
 Let $k\colon
\alpha \mapsto k_\alpha$ be a (non-negative) multiplicity function.
Assume the standard normalization of the $W$-invariant inner product.
We  put $u_\alpha = 1$ (type $A$), $u_\alpha = 
{\langle \alpha, \alpha \rangle \over 2}$ (type $B,C,D$)
for all $\alpha\in\Sigma$. Set $q_\alpha = q^{u_\alpha}$,
$t_\alpha = q^{k_\alpha}$ ($\alpha\in \Sigma$).
We define an inner product on functions on the torus $T$
by putting:
$$
\Delta^+ := \prod_{\alpha\in \Sigma^+} 
{(e^\alpha;q_\alpha)_\infty \over (t_\alpha e^\alpha; q_\alpha)_\infty},
\quad \Delta := \Delta^+ \overline{\Delta^+}, \quad
\langle f, g\rangle_k := \int_T f(x) \overline{g(x)} \Delta(x) dx.
\eqno\eq{Macdel}
$$
Using this inner product instead of \eqtag{clweight}, one now
proceeds in exactly the same way as in section 1 to
obtain (cf.\ \cite{\Ma}, \cite{\Mb}, \cite{\Mc})
{\sl Macdonald's symmetric polynomials} 
$P_\lambda^k\in \FSA(\Sigma)^W$ ($\lambda\in P^+$)
corresponding to the pair $(\Sigma, \Sigma)$ (type $A$) or
$(\Sigma, \Sigma^\vee)$ (type $B,C,D$).

In Koornwinder's case, it is convenient to use the notation
$x_i = e^{\eps_i}$, $x^\lambda = e^\lambda$ for $\lambda\in P$.
Recall that $P$ is the free $\ZZ$-span of 
the standard basis vectors $\eps_i\in P$. The definition of
Koornwinder's {\sl BC-type Askey-Wilson polynomials} is completely analogous
to Macdonald's case, except that one has to take a different
weight function $\Delta^+$ (cf.\ \eqtag{aworth}):
$$
\Delta^+(x) := \prod_{k=1}^l {(x_k^2;q)_\infty \over
(ax_k, bx_k, cx_k, dx_k;q)_\infty}\cdot
\prod_{i<j} {(x_i/x_j,x_ix_j;q)_\infty\over
(tx_i/x_j,tx_ix_j;q)_\infty} \; (a,b,c,d,t\in\CC).
\eqno\eq{KDel}
$$
For $n=1$ one reobtains the Askey-Wilson polynomials defined
in section 4.

To prove that the polynomials $P_\lambda$ are mutually 
orthogonal, Macdonald \cite{\Ma}, \cite{\Mb} and Koornwinder
\cite{\KWc} exhibit a self-adjoint partial 
$q$-difference operator $D_\sigma$
(depending on a minuscule weight $\sigma\in P^+$)
on $\FSA(\Sigma)^W$ which is diagonalized by the $P_\lambda$.
Here we take $\sigma=\eps_1$. The definition of
$D_\sigma$ then reads:
$$
\Phi_\sigma := {T_\sigma\Delta^+\over \Delta^+}, \quad
D_\sigma f :=
 |W\sigma|^{-1} \sum_{w\in W} (w\Phi_\sigma)(T_{w\sigma}f - f),
\quad f\in\FSA(\Sigma)^W.
$$
Here $T_\mu$ ($\mu\in P$) is the operator defined on $\FSA(\Sigma)$ by
$T_\mu x^\lambda = q^{\langle \lambda, \mu \rangle} x^\mu$
($\lambda\in P$).
The operator $D_\sigma$ maps $A(\Sigma)^W$ into itself and is diagonalized
by the $P_\lambda$. Koornwinder proves the orthogonality of the 
$P_\lambda$ under the assumption that (i) the parameters $a,b,c,d$ are real
or, if complex, appear in conjugate pairs, (ii) $|a|,|b|,|c|,|d| \leq 1$,
but the pairwise products of $a,b,c,d$ are not $\geq 1$, (iii) $t\in (-1,1)$.
In this case, the weight function $\Delta$ is continuous on the torus
$T$.  

Let us return to the setting of quantum symmetric spaces.
In all cases described above, Noumi and Sugitani have been able to
show that the radial part $D$ of a suitable Casimir operator
$C\in \FSZ\Uq$ (cf.\ \eqtag{Cdef}) 
is essentially the same as the partial $q$-difference
operator $D_\sigma$ corresponding to the restricted root system
$\Sigma=\Sigma(G,K)$. The multiplicity function
$k$ resp.\ the parameters $a,b,c,d,t$ are defined in terms 
of $q$ and the root multiplicities of $(G,K)$  in a way comparable
to \eqtag{multdef}. This result allows them to prove that the
zonal spherical functions $\phi(\lambda)$
($\lambda\in P^+_K$), when restricted  to the diagonal subgroup $\TT$,
can be expressed as 
Macdonald's polynomials or Koornwinder's Askey-Wilson polynomials.

The second case of the three cases alluded to above
arises when the embedding $\gok_\CC\subset \gog_\CC$
is induced by an embedding of Dynkin diagrams. 
By deleting one node in the Dynkin diagram of $\gog_\CC$, one 
obtains the (possibly disconnected) Dynkin diagram of a semisimple
complex Lie algebra $\gog_\CC'$. The Lie algebra $\gok_\CC := 
\gog_\CC'\oplus \CC$ naturally is a Lie subalgebra of $\gog_\CC$. 
The irreducible symmetric pairs $(G,K)$ 
that can be obtained via this procedure are 
exactly the irreducible Hermitian 
symmetric pairs\footnote{The author owes these observations to
Prof.\ Tom  H.\ Koornwinder.}. Observe that the Hermitian spaces
in Table I are those for which the root system is $BC_l$ and 
the long roots have multiplicity 1 (recall that $B_2\simeq C_2^\vee
\simeq C_2$ is
considered to be of type $BC$).

Since the
definition of $\Uq(\gog_\CC)$ is given in terms of generators
indexed by simple roots and relations depending only on
the Cartan matrix of $\gog_\CC$,
there will be a natural embedding of Hopf $\ast$-algebras
$\Uq(\gok_\CC)\subset \Uq(\gog_\CC)$ and a corresponding
surjective Hopf $\ast$-algebra morphism 
$\Aq(G) \to \Aq(K)$. From here, it is easy to define
a suitable coideal $\gok_q \subset \Uq(\gog_\CC)$ (cf.\ 
section 3). 
We shall call this
the {\sl regular embedded case}. It is quite opposite to
the maximally split case, since now the intersection
$\goh_\CC \cap \gok_\CC$ is maximal abelian in
$\gok_\CC$. For this reason, the restriction of 
$\gok_q$-biinvariant functions to the diagonal subgroup
$\TT$ will  not be injective. Therefore, the method 
proposed in section 3 does not work. 
To the author's knowledge, the spherical functions in
these cases have not yet been
studied, except in rank one (cf.\ \cite{\VS}, \cite{\Mas},
 \cite{\KWa}, \cite{\NYM}). 
It is likely, 
though, that the spherical functions in this case will
be some kind of polynomials orthogonal with respect to
a {\sl discrete} orthogonality measure. Until recently,
not much was known about multivariable orthogonal
polynomials with a discrete measure in connection with
root systems (however, see \cite{\St}).

The third case is a kind of {\sl interpolation} between 
the regular embedded and maximally split cases. It
arises when the Cartan algebra $\goh_\CC$ is not
necessarily invariant under the involution $\theta$, 
but its projection onto $\gop_\CC$ along the
decomposition $\gog_\CC = \gok_\CC \oplus \gop_\CC$ still
is maximal abelian in $\gop_\CC$. In this case, the restriction
of $\gok_q$-biinvariant functions to $\TT$ is still 
likely to be an injective operation. 

With this terminology, the coideal $\gok^\sigma$ defined
in section 4 is regular embedded ($\sigma = \pm \infty$),
maximally split ($\sigma = 0$), and interpolated 
($\sigma$ finite and non-zero). In the last case, 
the spherical functions are orthogonal with respect
to a mixed continuous and discrete orthogonality 
measure.

Very recently, Noumi, Sugitani and the author constructed
a one-parameter family of coideals $\gok^\sigma\subset
\Uq(\gog\gol(n,\CC))$ ($-\infty < \sigma < \infty$) defining
a quantum analogue of the symmetric space $U(n)/U(l)\times U(n-l)$
($l\leq [{n\over 2}]$).
 This symmetric space is different
from the others in the sense that the ``lowest'' spherical
representation occurs in the decomposition of the tensor
product $V\ten V^\ast$. In this case, the reflection
equation reads
$$
R^+ J_2 R^-_{21}J_1 = J_1 R^- J_2 R^+_{21},
\eqno\eq{Urefl}
$$
with $R^\pm_{21}:= PR^\pm P$. The coideal 
$\gok^\sigma$ is by definition spanned by the 
coefficients of the matrix $M:= L^+J^\sigma - J^\sigma L^-$,
where
$$
J^\sigma := \sum_{1\leq k\leq l}
 q^\sigma (q^{-\sigma} -q^\sigma)e_{kk}+
\sum_{l<k<l'} e_{kk} - \sum_{k\leq l\ \text{or}\ k\geq l'}
q^\sigma e_{kk'},
$$
and $k' := n+1-k$ ($1\leq k\leq n$).
This matrix $J^\sigma$ is a solution of the reflection equation
\eqtag{Urefl}.  For $l=1$ we reobtain \eqtag{Jsigmadef}.
The radial part of the Casimir operator \eqtag{Cdef} in this case
is essentially the same as Koornwinder's $q$-difference
operator $D_{\eps_1}$, the parameters
$a,b,c,d,t$ depending on two continuous parameters
$\sigma$, $\tau$, and two discrete parameters
$n$, $l$ (besides $q$). For details see the forthcoming paper
\cite{\DNS}.


\Refs\widestnumber\key{NYM}
\ref\key\AA \by G.E. Andrews, R. Askey \paper Classical orthogonal
polynomials \inbook in: ``Polyn\^omes Orthogonaux et Applications''
\bookinfo ed. C. Brzezinski, A. Draux, P. Magnus, P. Maroni, A. Ronveaux,
Lecture Notes in Mathematics 1171 \publ Springer-Verlag
\yr 1985\pages 36-62\endref
\ref\key\AW \by R.  Askey, J.  Wilson
\paper Some basic hypergeometric
orthogonal polynomials that generalize Jacobi polynomials
\jour Mem. Amer. Math. Soc.  \vol 54 \yr 1985 \issue 319\endref
\ref\key\Chin \by W. Chin, I.M. Musson \paper The coradical filtration
for quantized universal enveloping algebras \paperinfo preprint\endref
\ref\key\DKa \manyby M.S.  Dijkhuizen, T.H.  Koornwinder
\paper Quantum homogeneous spaces, duality 
and quantum 2-spheres
\jour Geom.  Dedicata \vol 52\yr 1994
\pages 291-315\endref
\ref\key\DKb \bysame
\paper CQG algebras: a direct algebraic approach
to compact quantum groups\jour Lett.  Math. Phys. \vol 32
\yr 1994 \pages 315-330\endref
\ref\key\DN \by M.S. Dijkhuizen, M. Noumi
 \paper A family of quantum projective 
spaces and related $q$-hyper\-geo\-metric orthogonal polynomials
\paperinfo preprint (1995)\endref
\ref\key\Dr \by V.G. Drinfel'd  \paper Quantum groups
\inbook in: Proceedings ICM Berkeley (1986) \bookinfo ed. A.M. Gleason
\pages 798-820 \publ Amer.  Math.  Soc. \publaddr Providence, RI \yr 1986
\endref
\ref\key\GR \by G. Gasper, M. Rahman
\book Basic hypergeometric series \bookinfo
Encyclopedia of Mathematics and 
Its Applications 35 \publ  Cambridge University Press
\publaddr London \yr 1990\endref
\ref\key\Hay \by T. Hayashi \paper Non-existence of homomorphisms between
quantum groups \jour Tokyo J. Math. \toappear\endref
\ref\key\HS \by G.J. Heckman, H. Schlichtkrull \book Harmonic 
Analysis and Special Functions on Symmetric Spaces \bookinfo 
Perspectives
in Math. 16\publ Academic Press \yr 1995\endref
\ref\key\Hela \manyby S. Helgason \book Differential Geometry, Lie Groups,
and Symmetric Spaces \bookinfo Pure and Applied Mathematics 80
\publ Academic Press \yr1978\endref
\ref\key\Helb \bysame \book Groups and Geometric Analysis
\bookinfo Pure and Applied Mathematics 113 \publ Academic Press \yr 1984
\endref
\ref\key\Jima \manyby M. Jimbo \paper A $q$-difference analogue of $U(\gog)$
and the Yang-Baxter equation\jour Lett. Math. Phys. 
\vol 10 \yr 1985 \pages 63-69 \endref
\ref\key\Jimb \bysame \paper Quantum $R$-matrix for the generalized
Toda system \jour Comm. Math. Phys. \vol 102 \yr 1986 \pages 537-547
\endref
\ref\key\Jimc \bysame \paper A $q$-analogue of $U(\gog\gol(n))$,
Hecke algebra and the Yang-Baxter equation
\jour Lett. Math. Phys.
\vol 11 \yr 1986 \pages 247-252\endref
\ref\key\KWa \manyby T.H. Koornwinder \paper Representations of the
twisted $SU(2)$ quantum group and some $q$-hypergeometric 
polynomials \jour Proc. Kon. Ned. Akad. Wet. Series A
\vol 92  \yr 1989 \pages 97-117 \endref
\ref\key\KWb \bysame 
\paper Orthogonal polynomials in connection 
with quantum groups \inbook in: ``Orthogonal polynomials: 
Theory and Practice'' \bookinfo ed. P. Nevai, NATO-ASI Series C 
294 \publ Kluwer \publaddr Norwell, MA \yr 1990\pages 257-292\endref
\ref\key\KWc \bysame \paper Askey-Wilson polynomials
for root systems of type $BC$ \inbook in: ``Hypergeometric functions on
domains of positivity, Jack polynomials, and applications''
\bookinfo ed. D.S.P. Richards 
Contemp. Math. 138 \publ Amer. Math. Soc. \publaddr 
Providence, RI \yr 1992 \pages 189-204\endref
\ref\key\KWd \bysame \paper Askey-Wilson polynomials
as zonal spherical functions on the $SU(2)$ quantum group
\jour SIAM J. Math. Anal. \vol 24 \yr 1993 \issue 3
\pages 795-813\endref
\ref\key\VK \by L.I. Korogodsky, L.L. Vaksman
\paper Quantum $G$-spaces and Heisenberg algebra
\inbook in: ``Quantum Groups'' \bookinfo 
ed. P.P. Kulish, Lecture Notes in Math. 1510
\publ Springer-Verlag \yr 1992 \pages 56-66\endref
\ref\key\Ku \by P.P. Kulish  \paper Quantum groups and quantum algebras
as symmetries of dynamical systems \paperinfo preprint YITP/K-959 \yr 1991
\endref
\ref\key\Lz \by G. Lusztig
\paper Quantum deformations of certain simple
modules over enveloping algebras \jour Adv. Math. \vol 70
\yr 1988 \pages 237-249 \endref
\ref\key\Ma \manyby I.G. Macdonald 
\paper A new class of symmetric functions
\inbook S\'eminaire Lotharingien de Combinatoire \bookinfo 
ed. L. Cerlienco, D. Foata, Publication I.R.M.A. 372/S-20,
Strasbourg \yr 1988 \pages 131-171\endref
\ref\key\Mb \bysame \paper Orthogonal polynomials 
associated with root systems \paperinfo preprint \yr 1988 \endref
\ref\key\Mc \bysame \paper Affine Hecke algebras and orthogonal
polynomials \jour S\'eminaire Bourbaki \vol 47\yr 1994-95
\issue 797\pages 1-18\endref
\ref\key\Mas \by T. Masuda, K. Mimachi, Y. Nakagami,
M. Noumi, K. Ueno \paper Representations of quantum groups and
a $q$-analogue of orthogonal polynomials \jour C. R. Acad. Sci. 
Paris S\'er. I Math. \vol 307 \yr 1988 \pages 559-564\endref
\ref\key\Nmac \by M. Noumi
\paper Macdonald's symmetric polynomials
as zonal spherical functions on some 
quantum homogeneous spaces
\jour Adv. Math. \toappear\endref
\ref\key\DNS \by M. Noumi, M.S. Dijkhuizen, T. Sugitani
\paper Multivariable Askey-Wilson polynomials and quantum complex
Grassmannians \paperinfo in:\ Proceedings of a Workshop on
$q$-Series, Special Functions and Related Topics, Toronto (1995)
 \jour Fields Inst. Comm. \toappear\endref
\ref\key\NSa \manyby M.  Noumi, T. Sugitani 
\paper Quantum symmetric spaces and related 
$q$-orthogonal polynomials \inbook in:
``Group Theoretical Methods in Physics''
\bookinfo Proceedings XX ICGTMP, 
Toyonaka (Japan), 1994, ed.\ A.\  Arima et al.
\publ World Scientific \publaddr Singapore
\yr 1995 \pages 28-40\endref
\ref\key\NSb \bysame \paperinfo in preparation\endref
\ref\key\NYM \by M.  Noumi, H.  Yamada, K.  Mimachi
\paper Finite-dimensional representations of the 
quantum group $GL_q(n,\CC)$
and the zonal spherical functions on 
$U_q(n-1)\backslash U_q(n)$
\jour Japanese J.  Math. \vol 19 \yr 1993 \issue 1 
\pages 31-80\endref
\ref\key\Po \by P. Podle\'s  \paper Quantum spheres
\jour Lett. Math. Phys. \vol 14 \yr 1987 
\pages 193-202 \endref
\ref\key\FRT \by N. Reshetikhin, L.D. Faddeev, 
L.A.  Takhtadjan  \paper Quantization of Lie groups 
and Lie algebras \jour Leningrad
Math. J. \vol 1 \yr 1990 \pages 193-225 \endref
\ref\key\Ro \by M.  Rosso
\paper Finite-dimensional representations of
the quantum analog of a complex simple Lie algebra
\jour Comm. Math. Phys.
\vol 117 \yr 1988 \pages 581-593\endref
\ref\key\St \by J. Stokman \paper Multivariable
big and little $q$-Jacobi polynomials\paperinfo Mathematical
Preprint Series 95-16,  University of Amsterdam\yr 1995\endref
\ref\key\Su \by T. Sugitani \paperinfo  in preparation\endref
\ref\key\UT \by K. Ueno, T. Takebayashi \paper Zonal 
spherical functions on quantum symmetric spaces and
Macdonald's symmetric polynomials
\inbook in: ``Quantum Groups'' \bookinfo  
ed. P.P. Kulish, Lecture Notes in Math. 1510
\publ Springer-Verlag \yr 1992 \pages 142-147 \endref
\ref\key\VS \by L.L. Vaksman, Y.S. So\u\i bel'man
\paper Algebra of functions
on the quantum group $SU(2)$ \jour Funct. Anal. Appl.
\vol 22 \yr 1988 \pages 170-181\endref
\ref\key\Wo \by S.L. Woronowicz \paper Compact matrix pseudogroups
\jour Comm. Math. Phys. \vol 111 \yr 1987 \pages 613-665\endref
\endRefs
\enddocument